\documentclass{elsart}
\usepackage{amssymb,amsmath}
\journal{The Journal of Algebra}

\theoremstyle{definition}
\newtheorem{BC}{Bounded Condition}

\begin{document}
\begin{frontmatter}

\title{Grothendieck rings of o-minimal expansions of ordered abelian groups}
\author[Kageyama]{M. Kageyama\corauthref{cor}},
\corauth[cor]{Corresponding author.}
\ead{kageyama@math.kyoto-u.ac.jp}
\author[Fujita]{M. Fujita}
\ead{fujita@math.kyoto-u.ac.jp}

\address[Kageyama]{Department of Mathematics,
Graduate School of Science,
Kyoto University,
Sakyou, Kyoto, 606-8502, Japan}
\address[Fujita]{Department of Mathematics,
Graduate School of Science,
Kyoto University,
Sakyou, Kyoto, 606-8502, Japan}

\begin{abstract}
We will calculate completely the Grothendieck rings, 
in the sense of first order logic, of o-minimal expansions of ordered abelian groups
by introducing the notion of the bounded Euler characteristic.
\end{abstract}
\end{frontmatter}

%\usepackage[all]{xy}

%----------------------------------- macros

%-----------------------------------

\section{Introduction}

The notion of the Grothendieck ring for a first-order structure was 
introduced by \cite{Kra-Thomas} and \cite{DL2}, independently.
In \cite{Kra-Thomas}, J. Kraj\'{\i}\v{c}ek and T. Scanlon 
clarified the relation between the triviality of this ring
and the non-existence of nontrivial weak Euler characteristic maps.
More precisely, they used weak Euler characteristics and Grothendieck rings
to handle the following situations. For instance, for a finite model and
when any one-to-one function is onto (PHP, pigeonhole principle),
however, for an infinite model, this dose not holds in general.
In \cite{DL2}, J. Denef and F. Loeser showed that for $T$ the theory of 
algebraically closed field containing a fixed field $k$,
it coincides  with the notion of the 
Grothendieck ring of algebraic varieties over $k$.
They treated with the motivic integration 
which was introduced by M. Kontsevich. 

For an arbitrary $\mathcal{L}$-structure $\mathcal{M}$,
$K_{0}({\mathcal{M}})$
and $K_{0}(M,{\mathcal{L}})$ denote the Grothendieck ring of the
$\mathcal{L}$-structure $\mathcal{M}$.

In \cite{vdD}, \cite{Clu-Hask} and \cite{Clu},
the Grothendieck rings of fields are calculated explicitly as follows:

\begin{enumerate}
\item $K_0(R, {\mathcal{L}}_{or})=\mathbb{Z}$, where
$R$ is a real closed field and ${\mathcal{L}}_{or}$ 
is the language $(<, +, -, \cdot, 0, 1)$.
\item $K_0(\mathbb{Q}_p, {\mathcal{L}}_{ring})=0$, where $p$ is a prime number, 
$\mathbb{Q}_p$ is the $p$-adic number field and
${\mathcal{L}}_{ring}$ is the language $(+, -, \cdot, 0, 1)$.
\item $K_0(\mathbb{F}_p((t)), {\mathcal{L}}_{ring})=0$, where $p$ is a prime number and
$\mathbb{F}_p((t))$ is the quotient field of the formal power series 
in the indeterminate $t$ over the finite field $\mathbb{F}_p$.
\item $K_0(F, {\mathcal{L}}_{ring})=0$, where $F$ denotes Laurent series fields $L((t_1))$,
$L((t_1))((t_2))$, $L((t_1))((t_2))((t_3))$ and $L$ is a finite extension of $\mathbb Q_p$ or $\mathbb F_q$.
Here $p$ is a prime number and  $q$ is a power of $p$.
\end{enumerate}

In \cite{Kra-Thomas} and \cite{DL2}, it is shown that the Grothendieck ring
$K_{0}(\mathbb{C}, {\mathcal{L}}_{ring})$ is extremely big and complicated:
\begin{enumerate}
\item[(5)] There exists a ring embedding
$\mathbb{Z}[X_j\;|\; j\in \mathfrak{c}] \hookrightarrow 
K_{0}(\mathbb{C}, {\mathcal{L}}_{ring})$,
where $\mathfrak{c}$ is the cardinality of continuum and 
$X_j$ $(j \in \mathfrak{c})$ are indeterminates.
\end{enumerate}

Although the Grothendieck rings of some structures have been calculated
as above, many other Grothendieck rings are not known yet
and the Grothendieck rings of o-minimal expansions of ordered abelian groups are known only a little.
See \cite{vdD} for the precise definition of an o-minimal structure.

In the present paper, we will calculate the Grothendieck rings of o-minimal expansions of ordered abelian groups
completely, namely, we have the following theorem:

\begin{thm}
\label{maintheorem}
Let ${\mathcal{G}}=(G, <, +, 0,\ldots)$ be an o-minimal expansion of an ordered abelian group.
Then $K_0(\mathcal{G})$ is isomorphic to either 
$\mathbb{Z}$ or the quotient ring $\mathbb{Z}[T]/(T^2+T)$ 
as a ring, where
$\mathbb{Z}[T]$ is a polynomial ring in an indeterminate $T$ over 
$\mathbb{Z}$ and $(T^2+T)$ is the ideal of $\mathbb{Z}[T]$ generated
by $T^2+T$.
\end{thm}

\section{Grothendieck Rings}

Let $\mathcal{M}$ be an $\mathcal{L}$-structure.
The notation ${Def^{n}(\mathcal{M})}$ denotes the family of all definable subsets
of $M^n$. We set
$\displaystyle {Def(\mathcal{M})}:=\bigcup_{n=0}^{\infty}{Def^{n}(\mathcal{M})}$.
Two definable sets $A, B \in {Def(\mathcal{M})}$ are {\it definably isomorphic},
denoted by $A \cong B$, if there is a definable bijection $A \rightarrow B$.

\begin{defn}[Grothendieck Ring]
\rm
\label{def-Grothendieck}
The {\it Grothendieck group} of an $\mathcal{L}$-structure $\mathcal{M}$ is
the abelian group $K_0(\mathcal{M})$ generated by symbols $[X]$, where
$X \in{Def(\mathcal{M})}$ with the relations $[X]=[Y]$ if $X$ and $Y$ are
definably isomorphic, and $[U\cup V]=[U]+[V]$ where $U$, $V\in {Def^{n}(\mathcal{M})}$ and
$U\cap V=\varnothing$. The ring structure is defined by
$[X][Y]=[X\times Y]$ where $X\times Y$ is the Cartesian product of definable
sets. The ring $K_0(\mathcal{M})$ with this multiplication is called 
{\it Grothendieck ring} of the $\mathcal{L}$-structure $\mathcal{M}$.
\end{defn}

\begin{rem}
\rm
\label{universal}
By construction, the map 
$[\;]:{Def(\mathcal{M})} \rightarrow K_0(\mathcal{M})$ satisfies
the following universal mapping property: \\
Consider the map $\chi:{Def(\mathcal{M})} \rightarrow \mathbb{Z}$ with
\begin{enumerate}
\item $\chi(U\cup V)=\chi(U)+\chi(V)$
for $U,V\in {Def^{n}(\mathcal{M})}$ with $U\cap V=\varnothing$, 
\item $\chi(X\times Y)=\chi(X)\cdot\chi(Y)$
for $X, Y\in {Def(\mathcal{M})}$,
\item $\chi(Z)=\chi(Z')$ if $Z, Z'\in {Def(\mathcal{M})}$, $Z\cong Z'$.
\end{enumerate}
Then, there exists an unique ring homomorphism 
$\psi :K_0(\mathcal{M})\rightarrow \mathbb{Z}$ such that 
$\psi\circ [\;] =\chi$.
\end{rem}

\begin{rem}
\rm
The onto-pigeonhole principle {\it ontoPHP} is the statement that
there is no set $A$, $a\in A$ and injective map $f$ from $A$ onto
$A\backslash \{a\}$.
By the construction of the Grothendieck ring of a structure $\mathcal{M}$,
$K_0(\mathcal{M})$ is nontrivial
if and only if ${\mathcal{M}}\models ontoPHP$.
See \cite{Kra-Thomas} for the details.
\end{rem}

\section{Grothendieck Rings of O-minimal Expansions of Ordered Abelian Groups}

We begin with the introduction of notations of 
an o-minimal structure $(G, <,\ldots)$. 

For a definable set $X\subseteq G^m$, we put
\begin{eqnarray*}
{\mathcal{C}}(X) &:=& \{f:X \rightarrow G \;|\; f \mbox{ is definable and continuous}\}, \\{\mathcal{C}}_{\infty}(X) &:=& {\mathcal{C}}(X)\cup\{-\infty,+\infty\},
\end{eqnarray*}
where we regard $-\infty$ and $+\infty$ as constant functions on $X$.
For $f\in \mathcal{C}(X)$, the graph of $f$ is denoted by
$\Gamma(f)\subseteq X\times G$.
For $f, g \in {\mathcal{C}}_{\infty}(X)$, we write $f<g$ if $f(x)<g(x)$ for all $x\in X$, 
and in this case we put
\[
(f, g)_X :=\{(x, r) \in X \times G \;|\; f(x) < r < g(x)\}.
\]

We next show that the Grothendieck rings of 
o-minimal expansions of ordered abelian groups are of the simple form:

\begin{lem}
\label{lemma1}
Let $(G, <, +, 0,\ldots)$ be an o-minimal expansion of an ordered abelian group. Then, 
\[
K_0(G)=\mathbb{Z}[\;[C]\;|\; C \subseteq G\;
\mbox{\rm is a cell}].
\]
\end{lem}

\begin{pf}
Let $M\subseteq G^n$ be a definable set.
By the cell decomposition theorem, 
\[M=C_1\cup\cdots\cup C_l\]
where $C_1,\ldots, C_l$ are cells. Hence
\[[M]=[C_1]+\cdots+[C_l].\]
Therefore, it suffices to show that for every cell $C \subseteq G^n,$
$[C]\in \mathbb{Z}[\;[C]\;|\; C \subseteq G\;
\mbox{\rm is a cell}]$. We will prove this by induction on $n$.
For simplicity we denote ${\mathbb{Z}}_{\rm cell}
:=\mathbb{Z}[\;[C]\;|\; C \subseteq G\;\mbox{\rm is a cell}]$.

The claim obviously holds true in the case where $n=1$.
Assume that the claim is true for $n=k$,
and we show that it holds for $n=k+1$.
Let $C\subseteq G^{k+1}$ be a cell.

If
\[C=\{(x,t) \in A \times G \;|\; t=f(x)\}\]
where $A\in G^k$ is the image $\pi(C)$ of $C$ under the projection 
$\pi:G^{k+1} \rightarrow G^k$ on the first $k$-coordinates and for some 
function $f\in {\mathcal{C}}(A)$. 
Hence there exist a definable bijection $C\cong A$.
Because $A$ is a cell, by the inductive assumption, 
$[C]=[A]\in {\mathbb{Z}}_{\rm cell}.$

If
\[C=\{(x,t) \in A \times G \;|\; \alpha(x)< t <\beta(x)\}\]
where $A\in G^k$ is the image $\pi(C)$ of $C$ under the projection
$\pi:G^{k+1} \rightarrow G^k$ on the first $k$-coordinates
and for some functions $\alpha, \beta\in {\mathcal{C}}_{\infty}(A)$.
\smallskip

{\bf{Case 1}:}
$\alpha=-\infty$, $\beta=+\infty$.

Then $C=A\times(-\infty, +\infty)$. Hence
$[C]=[A]\cdot[(-\infty,+\infty)]\in{\mathbb{Z}}_{\rm cell}.$
\smallskip

{\bf{Case 2}:}
$\alpha \in {\mathcal{C}}(A)$, $\beta=+\infty$.

Then we have a definable bijection,
\begin{center}
\begin{tabular}{ccc}
$A \times (0,+\infty)$&$\longrightarrow$&$C$ \\[-2.0mm]
$$&&$$ \\[-2.0mm]
$(x,t)$&$\longmapsto$&$(x,\alpha(x)+t)$.
\end{tabular}
\end{center}
Hence, 
$[C]=[A]\cdot[(0,+\infty)]\in{\mathbb{Z}}_{\rm cell}.$
\smallskip

{\bf{Case 3}:}
$\alpha =-\infty$, $\beta\in {\mathcal{C}}(A)$.

Then we have a definable bijection,
\begin{center}
\begin{tabular}{ccc}
$A \times (0,+\infty)$&$\longrightarrow$&$C$ \\[-2.0mm]
$$&&$$ \\[-2.0mm]
$(x,t)$&$\longmapsto$&$(x,\beta(x)-t)$.
\end{tabular}
\end{center}
Hence, 
$[C]=[A]\cdot[(0,+\infty)]\in{\mathbb{Z}}_{\rm cell}.$
\smallskip

{\bf{Case 4}:}
$\alpha$, $\beta\in {\mathcal{C}}(A)$.

Then,
\[C\cup\Gamma(\alpha)\cup D
=\{(x,t)\in A\times G \;|\; t< \beta(x)\}\]
where $D=\{(x,t)\in A\times G \;|\; t< \alpha(x)\}.$
Hence, by considering the case 3
\[[C]+[\Gamma(\alpha)]+[D]
\in {\mathbb{Z}}_{\rm cell}.\]
Because $[\Gamma(\alpha)],[D]\in {\mathbb{Z}}_{\rm cell},$ 
thus $[C] \in {\mathbb{Z}}_{\rm cell}.$ $\Box$
\end{pf}

\begin{cor}
\label{elementary-coro}
Let $(G, <, +, 0,\ldots)$ be an o-minimal expansion of an ordered abelian group.
We set $X:=[(0,+\infty)]$.
Then the equation $X^2+X=0$ holds true, and
\[K_0(G)=\{m+nX \;|\; m,n \in \mathbb{Z}\}.\]
\end{cor}

\begin{pf}
First, we prove the following claim.

\begin{claim}
\begin{enumerate}
\item[$i.)$] For the interval $(a, b)$ where $a, b \in G$,
$[(a, b)]=-1$,
\item[$ii.)$] $[(-\infty, +\infty)]=2X+1$,
\item[$iii.)$] $X^2=-X$.
\end{enumerate}
\end{claim}

{\it Proof of Claim.}

i.)\ 
Because $(a, b)\cong (0,b-a)$, 
we may assume $a=0$ and show that $[(0, b)]=-1$.
$(0, b) \cong (0, b/2) \cong (b/2, b)$
and, $[(0,b)]=[(0,b/2)]+1+[(b/2,b)]$. Hence, $[(0,b)]=-1$.
\smallskip

ii.)\ 
$(-\infty,0) \cong (0,+\infty)$ and
$[(-\infty,+\infty)]=[(-\infty,0)]+1+[(0,+\infty)]$
thus $[(-\infty, +\infty)]=2X+1$.
\smallskip

iii.)\ 
Let $I$ be the interval $(0, +\infty)$ and
$f:I \rightarrow I (x \mapsto x)$ be a function.
Then, $\displaystyle I\times I=(0, f)_I \cup\Gamma(f)\cup (f, +\infty)_I$.
We can construct the following definable bijections,
\begin{center}
\begin{tabular}{ccc}
$(0,f)_I$&$\longrightarrow$&$(f,+\infty)_I$ \\[-2.0mm]
$$&&$$ \\[-2.0mm]
$(x,y)$&$\longmapsto$&$(y,x)$
\end{tabular}
and
\begin{tabular}{ccc}
$I \times I$&$\longrightarrow$&$(f,+\infty)_I$ \\[-2.0mm]
$$&&$$ \\[-2.0mm]
$(x,y)$&$\longmapsto$&$(x,x+y)$.
\end{tabular}
\end{center}
Because $\Gamma(f)\cong I$,
\begin{eqnarray*}
[I\times I]&=&
[(0,f)_I]+[\Gamma(f)]+[(f,+\infty)_I]\\
&=&[I\times I]+[I]+[I\times I].
\end{eqnarray*}
We get $[I\times I]+[I]=0$. Thus $X^2+X=0$. $\Box$

By Lemma \ref{lemma1}, 
for each element $F \in K_0(G)$ there exist cells $C_1,\ldots,C_n$ in $G$
such that 
$$F=\sum_{j_1,\ldots,j_n}a_{j_1,\ldots,j_n}[C_1]^{j_1}\cdots[C_n]^{j_n}$$
where $a_{j_1,\ldots,j_n} \in \mathbb{Z}$.
Each cell $C_i(i=1,\ldots, n)$ is a point or an interval
and $(0,+\infty)\cong (a,+\infty)\cong (-\infty,b)\cong(-\infty,0)$
where $a, b \in G$. Using the above claim, 
we obtain $F=m+nX$ for some $m, n\in \mathbb{Z}$. $\Box$
\end{pf}
\smallskip

Next we will define a class of definable sets for every o-minimal expansion of an ordered abelian group 
and show its useful properties to calculate the Grothendieck rings
of o-minimal expansions of ordered abelian groups.

\begin{defn}
\rm
Let $(G,<,+,0,\ldots)$ be an o-minimal expansion of an ordered abelian group.
We call that a definable set $M \subseteq G^n$ is
{\it bounded} if $M \subseteq [b,b']^n$ for some $b, b' \in G$,
where $[b,b']:=\{t\in G\;|\;b\leq t \leq b'\}$.
\end{defn}

\begin{lem}
\label{bounded D}
Let $(G, <, +, 0,\ldots)$ be an o-minimal expansion of an ordered abelian group and 
$M \subseteq G^n$ be a bounded definable set with $\dim{M}=1$.
Then, there exists a definable bijection $M\rightarrow D$ for some
bounded definable set $D \subseteq G$.
\end{lem}

\begin{pf}
Since $\dim{M}=1$, by the cell decomposition theorem
we get the following decomposition
\[M=C_1 \cup\cdots\cup C_l \cup C_{l+1} \cup\cdots\cup C_m\]
where $C_1,\ldots, C_m$ are cells, $\dim{C_{1}}=1,\ldots,\dim{C_{l}}=1$
and $\dim{C_{l+1}}=0,\ldots,\dim{C_{m}}=0$.

\begin{claim}
For all $i=1,\ldots,l$, there exists a projection
$p_{n_i}:G^n \rightarrow G ((x_1,\ldots,x_n)  \mapsto x_{n_i})$
for some $1 \leq n_i \leq n$ such that 
$p_{n_i}|C_i:C_i\rightarrow p_{n_i}(C_i)$ is definably bijective.
\end{claim}

{\it Proof of Claim.}

We prove this claim by the induction on $n$.
When $n=1$, because each $C_i$ is an interval or a point, 
the claim holds true. Under the assumption that the claim holds true for $n=k$,
we show that the claim holds for $n=k+1$.
Let $p_1:G^{k+1}\rightarrow G$ be the projection on the first 
coordinate.%$(x_1,\ldots,x_{k+1})\mapsto x_1$.
\smallskip

{\bf{Case 1}:}
$\dim{p_1(C_i)}=1$.

For the projections $\pi_q:G^{k+1} \rightarrow G^q (q=1,\ldots,k+1)$ 
on the first $q$-coordinates, $\dim{\pi_{q}}(C_i)=1$, because
$\dim{C_i} \geq \dim{\pi_q(C_i)} \geq \dim{p_1(C_i)}=1$.
Hence, each cell $\pi_q(C_i)\;(q=2,\ldots,k+1)$ is the graph of
a definable function $f_q \in {\mathcal{C}}(\pi_{q-1}(C_i))$.
By using $f_2,\ldots,f_k$, we inductively define functions
$g_2,\ldots,g_{k+1}:p_1(C_i) \rightarrow G$ as follows. 
$g_2(x):=f_2(x)$ and we define
$g_{j+1}$ by $g_{j+1}(x):=f_{j+1}(x,g_2(x),\ldots,g_j(x))$
where $2 \leq j \leq k+1$ and $x\in p_1(C_i)$.
Then, for a definable function $g:p_1(C_i)\rightarrow G^k\;
(x\mapsto (g_2(x),\ldots,g_{k+1}(x)))$, $C_i=\Gamma(g)$.
Thus we obtain a definable bijection $p_1|C_i:C_i\rightarrow p_1(C_i)$.
\smallskip

{\bf{Case 2}:}
$\dim{p_1(C_i)}=0$.

Since $\dim{p_1(C_i)}=0$,
there are a point $a_i\in G$ and a cell $D_i\subseteq G^k$ such that
$C_i=\{a_i\}\times D_i$.
By inductive assumption, there is a projection 
$p_{n_i}:G^{k}\rightarrow G$ such that $p_{n_i}|D_i$ is injective.
Let $\tau$ be a projection such that 
$\tau:G^{k+1} \rightarrow G^k ((x_1,\ldots,x_{k+1})
\mapsto (x_2,\ldots,x_{k+1}))$. Then, $p_{n_i+1}=p_{n_i}\circ\tau$
and $p_{n_i+1}|C_i :C_i\rightarrow p_{n_i}(C_i)$
is a definable bijection. $\Box$
\smallskip

By Claim, each $C_i\;(i=1,\ldots,l)$ is definably 
bijective to an interval of $G$ and each $C_i\;(i=l+1,\ldots,m)$
is a point set. Thus, we can define a definable bijection $M \rightarrow D$
for some bounded definable set $D\subseteq G$. $\Box$
\end{pf}

\begin{prop}
\label{interval-prop}
Let $\mathcal{G}=(G, <, +, 0,\ldots)$ be an o-minimal expansion of an ordered abelian group,
$M \subseteq G^m$ be a non-bounded definable set and
$N \subseteq G^n$ be a bounded definable set. 
If $M$ and $N$ are definably isomorphic, 
then there exists a definable bijection 
$(0, +\infty)\rightarrow D$
for some bounded definable set $D \subseteq G$.
\end{prop}

\begin{pf}
Let $\pi_q:G^n \rightarrow G^q$ be the projection on the first 
$q$-coordinates. By the cell decomposition theorem,
\[M=C_1 \cup\cdots\cup C_m\] where $C_1,\ldots, C_m$ are cells.
Since $M$ is a non-bounded definable set,
we can choose a non-bounded cell $C_i$ for some $1\leq i\leq m$.
Because $C_i$ is non-bounded %Necessary by replacing a axis,
we may assume that $\pi_1(C_i)$ is a non-bounded interval $I$. %=(0,+\infty)$.

If $\pi_2(C_i)=\Gamma(f)$ for some $f\in {\mathcal{C}}(\pi_1(C_i))$, then 
we can define a definable injection $i_2:I\rightarrow \pi_2(C_i)$ by
$i_2(x):=(x,f(x))$.

If $\pi_2(C_i)=\{(x,y)\in I\times G\;|\;\alpha(x)<y<\beta(x)\}$
for some $\alpha,\beta \in {\mathcal{C}}_{\infty}(\pi_1(C_i))$, 
note that $G$ is a vector space over $\mathbb{Q}$
\cite[Chapter 1, Proposition 4.2]{vdD}, 
we can define a definable injection $i_2:I\rightarrow \pi_2(C_i)$ by
$$i_2(x):=\left\{
    \begin{array}{ll}
    (x,x) & \mbox{\rm if $\alpha=-\infty, \beta=+\infty$,} \\
    (x,\beta(x)-a) 
     & \mbox{\rm if $\alpha=-\infty, \beta\in {\mathcal{C}}(\pi_1(C_i))$,}\\
    (x,\alpha(x)+a) & \mbox{\rm if $\alpha \in {\mathcal{C}}(\pi_1(C_i)), \beta=+\infty$,} \\
    (x,(\alpha(x)+\beta(x))/2)
    & \mbox{\rm if $\alpha \in {\mathcal{C}}(\pi_1(C_i)), \beta\in {\mathcal{C}}(\pi_1(C_i))$,}
    \end{array}
    \right. $$
where $a$ is a positive element of $G$.

By continuing in the similarly way, we get a sequence of definable injections
\[I \stackrel{i_2}{\to} \pi_2(C_i) \stackrel{i_3}{\to} \cdots \stackrel{i_{n-1}}{\to}
\pi_{n-1}(C_i) \stackrel{i_n}{\to} C_i.\]
Let $\iota:I\rightarrow C_i$ be the composition of these definable injections.
Because $\dim{f(\iota(I))}=1$ by Lemma 
\ref{bounded D}, there is a bounded definable set $D \subseteq G$ such that
$f(\iota(I)) \cong D$. Thus we get a definable bijection between $I$ and $D$.
$\Box$
\end{pf}

It is easier to calculate the Grothendieck ring of the structure
$\mathcal{G}$ in the case where a non-bounded definable set
and a bounded definable set are definably isomorphic than 
in the other case.
To treat the latter case, we rewrite the condition as follow:

\begin{BC}
\rm
Let $(G, <, +, 0,\ldots)$ be an o-minimal expansion of an ordered abelian group, 
$M \subseteq G^m$ be a bounded definable set 
and $N \subseteq G^n$ be a definable set. 
If $M$ and $N$ are definably isomorphic, then $N$ is bounded.
\end{BC}

\begin{exmp}
\rm
Let $\mathcal{G}=(G,+,-,<,0)$ be the ordered divisible abelian group.
Then $\mathcal{G}$ satisfies Bounded Condition.
\end{exmp}

\begin{pf}
Suppose not. Then there are definable sets $X\subseteq G^m$, 
$Y\subseteq G^n$ such that
$X$ is non-bounded, $Y$ is bounded and $X\cong Y$. By Proposition \ref{interval-prop},
there is a definable bijection $f:(0,+\infty) \rightarrow D$ for some
bounded definable set $D\subseteq G$. Because $\mathcal{G}$ is o-minimal,
we may assume that $D$ is an interval $(a, b)$ for some $a,b \in G$.
By the monotonicity theorem \cite[Chapter 3, Theorem 1.2]{vdD},
there are points $a_1< \cdots <a_k$ in $(0, +\infty)$ such
that on each subinterval $(a_j, a_{j+1})$ with $a_0=0, a_{k+1}=+\infty$,
the function $f|(a_j,a_{j+1})$ is strictly monotone and continuous.
Since $(g:=)f|(a_k, +\infty):(a_k,+\infty)\rightarrow (a, b)$ is definable
and the ordered divisible abelian group admits quantifier elimination
\cite[Chapter 1, Corollary 7.8]{vdD}, the definable function $g$ is a polygonal
line. By dividing suitably $(a_k, +\infty)$ again,
we obtain points $a'_{k+1}< \cdots <a'_{n}$
in $(a_k, +\infty)$ with $a'_{k}=a_{k},a'_{n+1}=+\infty$,
and linear functions
$g_{k, k+1}:(a'_{k},a'_{k+1})\rightarrow(a,b),\ldots,
g_{n,n+1}:(a'_{n},a'_{n+1})\rightarrow (a,b)$ with
$g_{k,k+1},\ldots,g_{n,n+1}$ are strictly monotone.

There exist $m,m'\in\mathbb{Z}$ such that
$g_{n,n+1}(x)=mx+m'$, $m\neq 0$ where $x\in (a'_{n},a'_{n+1})$. When $m>0$ for
$x_0\in G$ with $(-m'+b)/m \ll x_0$, $g_{n,n+1}(x_0)>b$.
This contradicts to the fact that the target space of $g_{n,n+1}$
is $(a,b)$. We can also lead a contradiction when $m<0$ in the 
same way. $\Box$
\end{pf}

\begin{exmp}
\rm
Let $\mathcal{R}=(R,+,-,\cdot,<,0,1)$ be a real closed field.
Then $\mathcal{R}$ dose not satisfy Bounded Condition.
\end{exmp}

\begin{pf}
We can define a definable bijection 
$\phi:(0,1)\rightarrow(1,+\infty)$ by $\phi(x):=x/(1-x)$. $\Box$
\end{pf}

\section{Bounded Euler Characteristic}

We first recall the definition of the geometric Euler characteristic
\cite[Chapter 4]{vdD}.

\begin{defn}
\rm
Let $(G, <,\ldots)$ be an o-minimal structure
and $S$ be a definable subset of $G^m$. There exists a finite partition
$\mathcal{P}$ of $S$ into cells ${\mathcal{P}}=\{C_1, \ldots, C_l\}$ by 
the cell decomposition theorem. Then we define
the {\it geometric Euler characteristic} of the definable set $S$:
\[\chi_{g}(S):=\sum_{C \in{\mathcal{P}}}(-1)^{\dim{C}}\]
\end{defn}

This definition is seem to depend on the partition $\mathcal{P}$
of $S$. However, the definition dose not  depend on the choice of 
finite partitions. Moreover, it is known that 
$\chi_g$ is invariant under definable bijections
and satisfies the properties (1), (2) and (3) in Remark \ref{universal}.
See \cite[Chapter 4]{vdD} for the details.

\begin{lem}
\label{inj}
Let $\mathcal{G}=(G, <, +, 0,\ldots)$ be an o-minimal expansion of an ordered abelian group.
Consider the ring homomorphism $i:\mathbb{Z}\rightarrow K_0(G)$ given
by $i(1)=[\text{\rm one point}]$.
Then $i$ is injective.
\end{lem}

\begin{pf}
Consider the geometric Euler characteristic 
$\chi_g:Def(\mathcal{G})\rightarrow\mathbb{Z}$. By Remark \ref{universal}
there exists a ring homomorphism $\psi_g:K_0(\mathcal{G})\rightarrow\mathbb{Z}$
such that $\psi_g\circ[\;]=\chi_g$.
Fix $n\in\ker(i)$. We may assume that $n\geq 0$.
By the definition of $\chi_g$, 
\[n=\chi_g(\mbox{\rm $n$ points})=\psi_g\circ i(n)=0.\]
We have shown that $i$ is injective. $\Box$
\end{pf}

By Lemma \ref{inj}, we may consider naturally that $\mathbb{Z}$
is a subring of $K_0(\mathcal{G})$ for each o-minimal expansion of an ordered abelian group 
$\mathcal{G}$.

\begin{defn}
\rm
Let $(G, <, +, 0,\ldots)$ be an o-minimal expansion of an ordered abelian group,
$C\subseteq G^n$ be a cell and $p_k:G^n\rightarrow G^k$ be the
projection on the first $k$-coordinates. A cell $C$ is called 
{\it exceptional} if there exist $k\in \mathbb{N}$ and a cell
$A\subseteq G^{k-1}$ with $p_k(C)=A\times G$. A non-exceptional
cell $C$ is called {\it bad} if there exist $k\in \mathbb{N}$ and a cell
$A\subseteq G^{k-1}$ with
\[
p_k(C)=\{(x,t)\in A\times G\;|\; t < f(x)\}\;
      \text{or}\;\{(x,t)\in A\times G\;|\; f(x) < t \},
\]
where $f:A\rightarrow G$ is a definable function.
A {\it good} cell $C$ is a cell which is not neither exceptional nor bad.
\end{defn}

\begin{lem}
Let $(G, <, +, 0,\ldots)$ be an o-minimal expansion of an ordered abelian group,
$X \subseteq G^n$ be a definable set, $\mathcal{F}$ be a finite partition
of $X$ into cells any one of whose cell is not exceptional.
We put 
\[{\chi_b}(X):=\left\{
    \begin{array}{cl}\displaystyle
    \sum_{C\in {\mathcal{F}},C{\rm:good}}(-1)^{\dim C}
    & \mbox{\rm if ${\mathcal{F}}$ includes a good cell}, \vspace{1mm}\\
    0 & \mbox{\rm otherwise}.
    \end{array}
    \right.\]
Then ${\chi_b}(X)$ dose not depend on the choice of the finite partition $\mathcal{F}$.
\end{lem}

\begin{pf}
We can take such a finite partition ${\mathcal{F}}=\{C\}$ of $X$ by applying
the cell decomposition theorem to definable sets
$X$, $\{(x_1,\ldots,x_n)\in G^n\;|\; x_i>0\}$,
$\{(x_1,\ldots,x_n)\in G^n\;|\; x_i=0\}$,
and $\{(x_1,\ldots,x_n)\in G^n \;|\; x_i<0\},(i=1,\ldots,n).$
We set 
\[\chi_{b}^{\mathcal F}(X):=
\sum_{C \in \mathcal F, C \text{:good}}(-1)^{\dim(C)}. \]
Let ${\mathcal{G}}=\{D\}$ be another partition. 
Our purpose of this proof is to show 
$\chi_{b}^{\mathcal G}(X)=\chi_{b}^{\mathcal F}(X)$. 
Let $\mathcal H$ be a finer partition than $\mathcal F$ and $\mathcal G$.
If $\chi_{b}^{\mathcal F}(X)=\chi_{b}^{\mathcal H}(X)$ and $\chi_{b}^{\mathcal G}(X)
=\chi_{b}^{\mathcal H}(X)$, then $\chi_{b}^{\mathcal G}(X)=\chi_{b}^{\mathcal F}(X)$. 
Hence we may assume that $\mathcal G$ is a finer partition than $\mathcal F$.
We prove $\chi_{b}^{\mathcal G}(X)=\chi_{b}^{\mathcal F}(X)$ by the induction on $n$.
Remark that 
\[\chi_{b}^{\mathcal F}(X)=
\chi_{g}\Bigl{(}\bigcup_{C \in \mathcal F, C \text{:good}}C\Bigl{)}
=\sum_{C \in \mathcal F, C \text{:good}} (-1)^{\dim(C)}.\]
We have only to show that, for any bad cell $C$ of $\mathcal F$, 
\begin{center}
$\displaystyle\ \sum_{D \in \mathcal G, D \subseteq C,D
\text{:good}} (-1)^{\dim(D)}=0$.
\end{center}
We fix $C \in \mathcal F$ and 
set 
\[E:=\bigcup_{D \in \mathcal G,D \subseteq C,D \text{:good}} D.\]
Remark that 
\[\sum_{D \in \mathcal G,D \subseteq C,D \text{:good}}
(-1)^{\dim(D)}= \chi_g(E).\]

When $n=1$, $E=(a,b]$, $E=[a,b)$ or $E=\varnothing$ for some $a,b \in G$. 
Hence $\chi_g(E)=0$. 

We consider the case where $n>1$. 
Let $p$ be the projection on the first $(n-1)$-coordinates.
Then $p(C)$ is a non-exceptional cell. 
Let $\mathcal G'=
\{D'\}$ be the family of all good cells of the form: 
$D'=p(D)$ for some $D \in \mathcal G$.
Set $F:= \displaystyle \bigcup_{D' \in \mathcal G'} D'$.
Consider two cases.
\vspace{2mm}

$\bullet$ First consider the case where $C$ is of the form:
\[
\{(x,t) \in p(C) \times G\;|\; t=f(x)\}\;
 \text{or}\;\{(x,t) \in p(C) \times G\;|\; f(x)<t<g(x)\}, 
\]
where $f,g:p(C) \rightarrow G$ are definable functions.  
Remark that $\chi_g(F)=0$ by the inductive hypothesis. 
Since $E=\{(x,t) \in F \times G\;|\; t=f(x)\}$ or 
$E=\{(x,t) \in F \times G\;|\; f(x)<t<g(x)\}$, $\chi_g(E)=0$.
\vspace{2mm}

$\bullet$ Consider the other case, then there exist definable functions
$f<g$ on $D' \in \mathcal G'$ such that 
\begin{eqnarray*}
E \cap p^{-1}(D')&=&\{(x,t) \in D' \times G\;|\; f(x)<t \leq g(x)\}, \\
&&\{(x,t) \in D' \times G\;|\; f(x) \leq t < g(x)\}\;\text{or}\;
\varnothing.
\end{eqnarray*}
In each case, $\chi_g(E \cap p^{-1}(D'))=0$. 
Since $\displaystyle E= \bigcup_{D' \in \mathcal F'} (E \cap p^{-1}(D'))$,
$\chi_g(E)=0$. $\Box$
\end{pf}

\begin{lem}\label{union}
Let $(G, <, +, 0,\ldots)$ be an o-minimal expansion of an ordered abelian group.
Let $X$ and $Y$ be definable sets. 
Then $\chi_b(X \cup Y) + \chi_b(X \cap Y)
=\chi_b(X)+\chi_b(Y)$.
\end{lem}

\begin{pf}
This lemma follows from the definition of $\chi_b$ obviously.
$\Box$
\end{pf}

\begin{prop}\label{bunkou}
Let $(G, <, +, 0,\ldots)$ be an o-minimal expansion of an ordered abelian group,
$X \subseteq G^{m+n}$ be a definable subset,
$\mathcal D$ be a decomposition of $G^{m+n}$ partitioning $X$
and $\pi:G^{m+n} \rightarrow G^m$ 
be the projection on the first $m$-coordinates.
Assume that all cells are not exceptional.
Given a cell $A \in \pi(\mathcal D)$ there is a constant 
$e_A$ with $\chi_b(X \cap p^{-1}(a))=e_A$ 
and $\chi_b(X \cap p^{-1}(A))=\chi_b(A)e_A$.
\end{prop}

\begin{pf}
Fix $A \in \pi(\mathcal D)$.
For each cell $C$ of $\mathcal D$, 
$C \cap \pi^{-1}(a)=\varnothing$ if $\pi(C) \not=A$ and $a \in A$.
If $\pi(C)=A$, $C \cap \pi^{-1}(a)$ is a cell and its dimension
does not depend on the choice of $a \in A$. 
Moreover, if $C \cap p^{-1}(a)$ is good for some $a \in A$,
the same statement holds true for all $a \in A$.
Set $e_A= \chi_b(X \cap \pi^{-1}(a))$ for some $a \in A$.
Then $e_A$ satisfies the requirement of the first statement of this lemma.
It is also obvious that 
$\chi_b(X \cap p^{-1}(A))=
\chi_b(A)e_A$ by the definition of $\chi_b$. $\Box$
\end{pf}

\begin{cor}\label{bunkou1}
Let $(G, <, +, 0,\ldots)$ be an o-minimal expansion of an ordered abelian group
and $X\subseteq G^m$ and $Y\subseteq G^n$ be definable sets. 
Then $\chi_b(X \times Y)=\chi_b(X) \cdot \chi_b(Y)$.
\end{cor}

\begin{pf}
This corollary follows from Proposition \ref{bunkou}. $\Box$
\end{pf}

\begin{lem}\label{sp}
Let $(G, <, +, 0,\ldots)$ be an o-minimal expansion of an ordered abelian group.
Moreover, assume that $\mathcal{G}$ satisfies Bounded Condition.
Then a cell $C$ is good if and only if $C$ is bounded.
\end{lem}

\begin{pf}
It is obvious that a cell which is not good is not bounded.
Hence we have only to show that a good cell $C \subseteq G^n$ is bounded. 
We prove it by the induction on $n$. When $n=1$, it is obvious. 
Consider the case when $n>1$. 
Let $p:G^n\rightarrow G^{n-1}$ be the projection on the first
$(n-1)$-coordinates. The cell $p(C)$ is bounded by the inductive hypothesis.
Let $d \in G$ such that $p(C) \subseteq [-d,d]^{n-1}$.
Remark that $C$ is of the form:
\[
\{(x,t) \in p(C) \times G\;|\;t=f(x)\}
\;\text{or}\;\{(x,t) \in p(C) \times G\;|\;f(x)<t<g(x)\},
\]
where $f$ and $g$ are definable functions on $p(C)$.
There exists positive $d' \in G$ such that 
$-d'<f(x)<d'$ and $-d'<g(x)<d'$ for all $x \in p(C)$.
Set $d'':= \max\{d,d'\}$.
Then $C \subseteq [-d'',d'']^n$, namely, $C$ is bounded. $\Box$
\end{pf}

\begin{lem}\label{koyuki}
Let $(G, <, +, 0,\ldots)$ be an o-minimal expansion of an ordered abelian group
satisfying Bounded Condition.
Let $X \subseteq G^m$ be a definable set and
$\sigma$ be a permutation of $\{1, \ldots,m\}$. 
We define a definable function 
$\Psi_{\sigma}: G^{m} \rightarrow G^{m}$ by 
$\Psi_{\sigma}(x_1, \ldots, x_{m})=(x_{\sigma(1)}, \ldots, x_{\sigma(m)})$. 
Then $\chi_b(X)=\chi_b(\Psi_{\sigma}(X))$.
\end{lem}

\begin{pf}
Since the symmetric group on $\{1, \ldots, m\}$
is generated by the transpositions $(i,i+1)$, 
we may assume that $\sigma=(i,i+1)$. 
By \cite[Chapter 4, Proposition 2.13]{vdD},
there exists a cell decomposition $\mathcal D$ such that 
any cell is not exceptional and $\Psi_{\sigma}(C)$ are
also cells for all cells $C \in \mathcal D$.
Since a cell is good if and only if it is bounded by Lemma 
\ref{sp}, $\Psi_{\sigma}(C)$ is good if and only if so is $C$.
Hence, $\chi_b(X)=\chi_b(\Psi_{\sigma}(X))$ by the definition of $\chi_b$.
$\Box$
\end{pf}

We are now ready to state the invariance of $\chi_b$
under bijections definable in an o-minimal expansion of an ordered abelian group which satisfies 
Bounded Condition.

\begin{prop}\label{bij}
Let $(G, <, +, 0,\ldots)$ be an o-minimal expansion of an ordered abelian group
satisfying Bounded Condition.
Let $X \subseteq G^m$ be a definable set and 
$f:X \rightarrow G^n$ be an injective definable map.
Then $\chi_b(X)=\chi_b(f(X))$. 
\end{prop}

\begin{pf}
Consider the graph $\Gamma(f) \subseteq G^{m+n}$ and the definable set
$\Gamma'(f)=\{(f(x),x) \in G^n \times X \}$.
By Proposition \ref{bunkou}, $\chi_b(X)=\chi_b(\Gamma(f))$
and $\chi_b(f(X))=\chi_b(\Gamma'(f))$.
Because $\chi_b(\Gamma(f))=
\chi_b(\Gamma'(f))$ by Lemma \ref{koyuki}. 
We obtain the conclusion. $\Box$
\end{pf}

\begin{defn}
\rm
\label{bdd-Euler}
Let $(G, <, +, 0,\ldots)$ be an o-minimal expansion of an ordered abelian group
satisfying Bounded Condition.
For all definable sets $X\subseteq G^n$,
we call $\chi_b(X)$ the {\it bounded Euler characteristic} of $X$.
\end{defn}

\begin{rem}
\rm
The following theorem ensures that our definition of $\chi_b$ 
coincides with the notion of the bounded Euler characteristic in 
\cite{Schanuel}.
\end{rem}

\begin{thm}\label{bd}
Let $\mathcal(G, <, +, 0,\ldots)$ be an o-minimal expansion of an ordered abelian group and
$X \in G^n$ be a definable set.
Let $d:X \rightarrow [0, \infty)$ be a definable function such that
$d^{-1}(t)$ is bounded for any $t \geq 0$.
Set $X_d(t):=\{x \in X\;|\;d(x)\leq t\}$ for any $t \in G$.
Then there exists $\mu \in G$ with $\chi_g(X_d(t))=\chi_b(X)$ for
$t\geq \mu$.
\end{thm}

\begin{pf}
Consider the definable set $\Gamma'(d):=\{(t,x) \in G \times X\;|\;d(x)=t\}$.
Let $p$ be the projection of $\Gamma'(d)$ to the first factor. 
Apply the cell decomposition theorem to $\Gamma'(d)$.
Let $\Gamma'(d)=C_1\cup\cdots\cup C_k$ be the cell decomposition. 
We may assume that $C_1, \ldots, C_j$ are bounded and $C_{j+1}, \ldots, C_k$
are not bounded.

Since the fibres of $d$ are bounded, the cell $C_i$ is bounded 
if and only if $p(C_i)$ is bounded.
Hence there exists $\mu \in G$ such that $C_i \cap p^{-1}(\{t \in G\;|\;t>\mu\})
= \varnothing$ for all $i=1, \ldots,j$ and 
$C_i \cap p^{-1}(\{t \in G\;|\;t>\mu\})\neq\varnothing$ for all $i=j+1, \ldots,k$. 
It is easy to see that the definable sets $C_i \cap p^{-1}(\{s \in G\;|\;s>t\})$
are cells of dimension $\dim C_i$ for all $i=j+1, \ldots,k$. 
Hence we omit the proof of this fact.

Fix $t \geq \mu$. Then  
\begin{align*}
\chi_g(X_d(t)) &=\chi_g(X)-\chi_g(\{x \in X\;|\;d(x)>t\}) \\
               &=\chi_g(X)-\sum_{i=j+1}^k (-1)^{\dim(C_i \cap p^{-1}
               (\{s \in G\;|\;s>t\}))}\\
&=\chi_g(X)-\sum_{i=j+1}^k (-1)^{\dim(C_i)} \;\;(\text{by the above fact}) \\
&=\sum_{i=1}^j(-1)^{\dim(C_i)} \\
&=\chi_b(\Gamma'(d)) \;\;(\text{by the definition})\\
&=\chi_b(X) \;\;(\text{by Proposition \ref{bij}}). \;\;\Box
\end{align*}
\end{pf}

\section{Proof of Theorem \ref{maintheorem}}

We are now ready to prove Theorem \ref{maintheorem}.

%\begin{Theorem}
%\label{main-theorem}
%Let ${\mathcal{G}}=(G, <, +, 0,\ldots)$ be an o-minimal expansion of an ordered abelian group.
%Then $K_0({\mathcal{G}})$ is isomorphic to either 
%$\mathbb{Z}$ or the quotient ring $\mathbb{Z}[T]/(T^2+T)$ as a ring, where
%$\mathbb{Z}[T]$ is a polynomial ring in an indeterminate $T$ over 
%$\mathbb{Z}$ and $(T^2+T)$ is the ideal of $\mathbb{Z}[T]$ generated by
%$T^2+T$.
%\end{Theorem}

\begin{pf}

{\bf{Case 1}:}
There exists a definable bijection between a non-bounded definable set
and a bounded definable set.

Then by Proposition \ref{interval-prop}, 
we can take a definable bijection $(0, +\infty) \cong D$ for some 
bounded definable set $D \subseteq G$. 
Because $[(0, +\infty)]=[D]\in \mathbb{Z}$, 
the ring homomorphism $i:\mathbb{Z}\rightarrow K_0(\mathcal{G})$ given by
$i(1)=[\text{\rm one point}]$ is surjective.
By Lemma \ref{inj}, $i$ is injective.
Therefore $K_0({\mathcal{G}})$ is isomorphic to $\mathbb{Z}$ as a ring.
\smallskip

{\bf{Case 2}:}
There exist no definable bijections of 
non-bounded definable sets into bounded definable sets.

Then, because $\mathcal G$ satisfies Bounded Condition,
we can define the bounded Euler characteristic $\chi_b$.
By Corollary \ref{elementary-coro},
the following ring homomorphism is surjective:
\begin{center}
\begin{tabular}{ccc}
$\phi:\mathbb{Z}[T]/(T^2+T)$&$\longrightarrow$&$K_0({\mathcal{G}})$ \\[-2.0mm]
$$&&$$ \\[-2.0mm]
$1$&$\longmapsto$&$[\text{\rm one point}]$\\[-.5mm]
$T$&$\longmapsto$&$X$
\end{tabular}
\end{center}
where $X=[(0,+\infty)]$.

We show that this ring homomorphism is injective.
Fix $m+nX \in\ker(\phi)$ where $m,n\in \mathbb{Z}$.
Considering the universal mapping property of $(K_0({\mathcal{G}}),[\;])$
for the geometric Euler characteristic $\chi_g$,
there exists an unique ring homomorphism 
${\psi}_g:K_0({\mathcal{G}})\rightarrow \mathbb{Z}$ such that 
$\psi_g \circ [\;]=\chi_g$.
By the definition of $\psi_g$,
\[
\psi_g(m+nX)=m+n\psi_g(X)
            =m+n\chi_g((0,+\infty))
            =m-n.
\]
Thus we get $m=n$. Similarly
for the bounded Euler characteristic $\chi_b$
there exists an unique ring homomorphism 
${\psi}_b:K_0({\mathcal{G}})\rightarrow \mathbb{Z}$ such that 
$\psi_b \circ [\;]=\chi_b$.
By the definition of $\psi_b$,
\[
\psi_b(m+nX)=m+n\psi_b(X)
            =m+n\chi_b((0,+\infty))
            =m.
\]
Thus we get $m=n=0$. We have shown $\phi$ is injective.
$\Box$
\end{pf}

\bibliographystyle{elsart-num}
\bibliography{kageyama}
\end{document}